\documentclass[12pt]{amsart}

\usepackage{amsmath}
\usepackage{amssymb}  
\usepackage{latexsym} 
\usepackage{comment}
\usepackage{url}
\usepackage{fullpage,url,amssymb,amsmath,amsthm,amsfonts,mathrsfs}
\usepackage[usenames,dvipsnames]{color}
\usepackage[pagebackref = true, colorlinks = true, linkcolor = blue, citecolor = Green]{hyperref}
\usepackage[alphabetic,lite]{amsrefs}
\usepackage{enumitem}
\usepackage{amscd}   
\usepackage[all, cmtip]{xy} 
\usepackage{xfrac}

\DeclareFontEncoding{OT2}{}{} 

\usepackage[all]{xy}
\usepackage{fullpage}

\usepackage{color} 

\def\act#1#2%
  {\mathop{}%
   \mathopen{\vphantom{#2}}^{#1}%
   \kern-3\scriptspace%
   #2}

\newcommand{\Z}{{\mathbb Z}}

\newcommand{\F}{{\mathbb F}}

\newcommand{\calI}{{\mathcal I}}

\newcommand{\calM}{{\mathcal M}}

\DeclareMathOperator{\Fr}{Fr}

\DeclareMathOperator{\Spec}{Spec}

\newtheorem{Theorem}{Theorem}[section]

\newtheorem{Conjecture}[Theorem]{Conjecture}
\newtheorem{Question}[Theorem]{Question}

\newtheorem{Example}[Theorem]{Example}

\numberwithin{equation}{section}

\begin{document}

\title{Maps between curves and arithmetic obstructions}

\author{Andrew V. Sutherland}
\address{Department of Mathematics, Massachusetts Institute of Technology, 77 Massachusetts Avenue, Cambridge, Massachusetts 02139, United States of America}
\email{drew@math.mit.edu}
\urladdr{http://math.mit.edu/~drew}

\author{Jos\'e Felipe Voloch}
\address{School of Mathematics and Statistics, University of Canterbury, Private Bag 4800, Christchurch 8140, New Zealand}
\email{felipe.voloch@canterbury.ac.nz}
\urladdr{http://www.math.canterbury.ac.nz/\~{}f.voloch}

\begin{abstract}
Let $X$ and $Y$ be curves over a finite field.  In this article we explore methods to determine whether there is a rational map from $Y$ to $X$ by considering $L$-functions of certain covers of $X$ and $Y$ and propose a specific family of covers to address the special case of determining when $X$ and $Y$ are isomorphic.
We also discuss an application to factoring polynomials over finite fields.
\end{abstract}

\maketitle
%


\section{Introduction}

Given two algebraic curves $X,Y$ of genus at least two
over a finite field we would like to decide if there is a rational map from $Y$ to $X$.
Hess and M\"ohlmann \cites{Hess, Moehlmann} have given an algorithm to decide if such a map exists by
performing an optimised search for the map up to a known bound; if the search is unsuccessful then no such map exists.
This procedure works well when the map exists, but it may be very time consuming when it does not.
The purpose of this paper is to provide a way of deciding when there is no such map without performing an exhaustive search. We concentrate on the case of isomorphisms but briefly touch on the general case.
Our methods can also provide a short certificate of the non-existence of a rational map, a feature not available with existing algorithms.

A result of Poonen \cite{Poonen} extending an idea of Kayal
shows that, given a one-parameter family of curves over a finite
field with distinct $L$-polynomials for distinct values of the parameter and a suitable bound on the genus,
one can construct a deterministic polynomial time algorithm for factoring polynomials over that field. Our investigations
suggest some candidate families of such curves, but unfortunately we cannot prove that they work.

\section{The fundamental group}

Let $X/K$ be a smooth geometrically connected variety over a field $K$.
Let $G_K$ be the absolute Galois group of $K$ and $\bar{X}$ the base-change of $X$ to an algebraic closure
of $K$. We denote by $\pi_1(.)$ the algebraic fundamental group functor
on (geometrically pointed) schemes and we omit base-points from the
notation. We have the fundamental exact sequence
\begin{equation}
\label{fund}
1 \rightarrow \pi_1(\bar{X}) \rightarrow \pi_1(X) \rightarrow G_K \rightarrow 1.
\end{equation}
The map $p_X: \pi_1(X) \rightarrow G_K$ from the above sequence is obtained by functoriality from the
structural morphism $X \to \Spec K$. Grothendieck's anabelian program is to specify a class of varieties,
termed anabelian, for which the varieties and morphisms between them can be recovered from the corresponding
fundamental groups together with the corresponding maps~$p_X$ when the ground field is finitely generated over
$\mathbb{Q}$. There has been some work done over finite fields as well, although the anabelian program will not
work in the same way (the analogue of the section conjecture is false, for example). 

For the rest of the paper we restrict to the case where $K$ is a finite field. 
As usual, $\F_q$ is the field of $q$ elements and we denote
by $p$ its characteristic. Here is a positive result.

\begin{Theorem}
\label{moch}
(Mochizuki--Tamagawa) Let $X,Y$ be smooth projective curves of genus at least two
over a finite field $\F_q$. If
there is an isomorphism from $\pi_1(X)$ to $\pi_1(Y)$ inducing the identity on $G_{\F_q}$ via $p_X,p_Y$, then
$X$ is isomorphic to $Y$.
\end{Theorem}

The fundamental group is a mysterious object. What kind of information
can we extract from it? First of all, if $J_X$ denotes the Jacobian of
a curve $X$ we have the fundamental exact sequence \eqref{fund} for $J_X$ also,
and $\pi_1(\bar{J_X})$ is the abelianisation of $\pi_1(\bar{X})$; thus
the prime-to-$p$ part 
of $\pi_1(\bar{J_X})$ is the product of the Tate modules of $J_X$.
The fundamental exact sequence describes the Galois action on the
Tate module, so its description is equivalent to the $L$-function of
$X$, which we can compute by counting points on $X$ (over suitable extensions of~$\F_q$). But by the very nature of the fundamental
group, we can count points on covers as well. Since knowing
$J_X$ alone up to isogeny is not enough to recover $X$, we need to pass
to covers. According to J. Stix (personal communication) the proof of Theorem \ref{moch}
requires only solvable covers. The most natural covers come from the Hilbert class field tower.
Let $\Fr: J_X \to J_X$ denote the $\F_q$-Frobenius map. Define 
$H(X) := (I-\Fr)^*(X) \subset J_X$; it is an unramified abelian cover of $X$ with
Galois group $J_X(\F_q)$, well defined up to a twist that corresponds to a choice of
divisor of degree one embedding $X$ into $J_X$. 
Define $H_0(X):=X$, $H_1(X):= H(X)$ and
successively define $H_{n+1}(X) := H_n(H(X))$ for integers $n\ge 1$. These covers
can be computed from $\pi_1(X)$ but are perhaps more computationally accessible.

\begin{Conjecture}\label{main-conj}
Let $X,Y$ be smooth projective curves of genus at least two
over a finite field~$\F_q$. If, for each $n$, there are choices of twists
such that the $L$-function of $H_n(X)$ is equal to
the $L$-function of $H_n(Y)$ for all $n \ge 0$, then $X$ is isomorphic to $Y$.
\end{Conjecture}

There are exactly $8$ curves for which $H(X)=X$, equivalently, curves whose function fields
have class number one \cites{Shen, Stirpe}. We have verified that Conjecture~\ref{main-conj}
holds when $X$ and $Y$ are among this list of $8$ curves, and it therefore holds if either $H(X)=X$ or $H(Y)=Y$.
We may thus assume henceforth that $X$ and $Y$ have non-trivial Hilbert class field towers.

The basis for our heuristic is the following consequence of the usual calculation leading to the
birthday paradox. If $\calM$ and $\calI$ are finite sets with cardinalities $M$ and $I$ respectively,
then for $I\preceq M^2$ the probability that a random map $\calM\to \calI$ is non-injective is
bounded above zero, but for $I$ asymptotically larger than $M^2$ this probability decays rapidly to zero.
Explicitly, this probability is $\prod_{j=0}^{M-1} (1-j/I) \sim e^{-M(M-1)/2I}$. 

We first apply this to the set $\calM$ of isomorphism classes of curves of some fixed genus $g>1$ 
over a finite field $\F_q$ and the set $\calI$ of isogeny classes of abelian varieties of dimension~$g$
over a finite field~$\F_q$.
For fixed $g$ and large $q$ we have $M \sim q^{3g-3}$ and $I \sim q^{g(g+1)/4}$, hence it is reasonable
to expect that there will be distinct curves with isogenous Jacobians when $6g-6 \le g(g+1)/4$, that is,
$g \le 22$ (see \cite{Howe-DiPippo} and \cite{Shankar-Tsimerman} where this kind of question is discussed in more detail).
For larger genus one expects this to be very rare, but we note a result of Mestre \cite{Mestre} that
allows one to construct, for every $g>1$, pairs of genus $g$ curves with isogenous Jacobians
(over sufficiently large finite fields).

Provided $H(X)\ne X$ (as we now assume), passing from $X$ to $H(X)$ increases the genus and thus
makes it more likely that we can use isogeny invariants to distinguish non-isomorphic curves.
For large $q$ and $g$, the genus of $H(X)$ is much larger than the genus $g$ of $X$;
indeed, it is on the the order of $gq^g$.
However, the Jacobian of $H(X)$ is not an arbitrary abelian variety of this dimension; it decomposes
up to isogeny as a product of abelian varieties of smaller dimension. The precise shape of the
decomposition depends on the group structure of $J(\F_q)$, but given its huge dimension, barring any
additional constraints, its isogeny class falls into a very small set of possibilities.
Thus on probabilistic grounds it is still likely that the map $X \mapsto L(H(X), T)$ is injective
when $g$ and $q$ are large.
This makes it plausible that, up to a finite set of exceptions, Conjecture~\ref{main-conj} holds
even when we restrict to $n\le 1$. The set of exceptions is non-empty as the examples in the next section show.

In addition to $H(X)$, we can also consider the cover of $X$ obtained by pulling back via multiplication by $2$ on
the Jacobian (assuming the characteristic is not $2$). 
This gives a cover $X^{(2)}$ of $X$ of degree $2^{2g}$. In general, this cover is not a subcover of the $H_n(X)$ 
considered above.
Its Jacobian $J_{X^{(2)}}$ is not a random abelian variety of its dimension, since it
decomposes (after base change to the algebraic closure of the ground
field) up to isogeny into a product of the Jacobian $J_X$ of $X$ and $2^{2g}-1$ abelian varieties of
dimension $g-1$ (the Prym varieties of $X$). If we assume that the isogeny classes of these factors
are random, we are picking them out of a set of size $\sim q^{{(2^{2g}-1)}g(g-1)/4}$ which is much smaller than
if we regarded $J^{(2)}$ as random, but still very large.
On the other hand we should note that the 
construction of Mestre \cite{Mestre} mentioned above produces curves that not only have
isogenous Jacobians, but a few of their Prym varieties will also be isogenous.

\begin{Question}\label{prym-quest} 
Are there non-isomorphic curves $X,Y$ over $\F_q$ of genus at least two and $p \ne 2$ with 
$J_{X^{(2)}}, J_{Y^{(2)}}$ isogenous?
\end{Question}

For maps between curves of different genera, it is less clear what to expect. In particular, we do not have a result generalizing Theorem \ref{moch}. But one can consider the following:

\begin{Question}\label{dommap-quest}
Let $X,Y$ be smooth projective curves of genus at least two
over a finite field~$k$, with $H(X)\ne X$ and $H(Y)\ne Y$.
Suppose the $L$-function of $H_n(X)$ divides
the $L$-function of $H_n(Y)$ for all $n \ge 0$.
Does this imply the existence of a dominant map $Y \to X$?
\end{Question}

As shown by an example of Brendan Creutz, the answer to Question~\ref{dommap-quest} is no if we allow $H(X)=X$.
A generalization of his idea (which is the case $n=0$) is as follows. Start with $X$ such that $X(\F_q) = \emptyset$. Consider
the Jacobian $J$ of $H_n(X)$ and by slicing with suitable hypersurfaces, construct a smooth 
curve $D \subset J$ with $0 \in D$, hence in particular $D(\F_q) \ne \emptyset$. This $D$ cannot
map to $X$ (as $X(k) = \emptyset$) but $L(H_n(X),t) | L(D,t) | L(H_n(D),t)$ by construction.
So the $n$ in Question~\ref{dommap-quest} cannot be uniformly bounded.

\section{Certifying non-isomorphism}

If two curves can be distinguished by the $L$-polynomials of low degree covers
then a succinct certificate can be given in the form of a prime $\ell$ for
which the corresponding two $L$-polynomials are distinct modulo $\ell$, 
together with the calculation of these polynomials; note that we can assume $\ell=O(g\log q)$,
since otherwise the $L$-polynomials must coincide.
For fixed $g$ the Schoof-Pila algorithm can be used to determine $\ell$ and compute
the $L$-polynomials modulo $\ell$ with a running time that is polynomial in $\log q$,
but exponential in $g$.
When~$g$ is large relative to $\log p$, where $p$ is the characteristic of $\F_q$, one is better off using algorithms based on $p$-adic cohomology to
compute the $L$-polynomials over $\Z$ and then reduce modulo a suitable prime $\ell$.
The complexity of the $p$-adic approach is polynomial in~$g$ but exponential in $\log p$.
The most general algorithm of this type is due to Tuitman \cite{Tuitman},
and is applicable to all curves that admit a suitable lift to characteristic zero; its
complexity is quasi-linear in $p$ and polynomial in~$g$.
When $q=p$ is prime one can instead apply Harvey's result for arithmetic schemes \cite{Harvey}, which
improves the dependence on $p$ to $O(p^{1/2+o(1)})$.
At present there is no algorithm known with a running time that is polynomial in
both $g$ and $\log q$, thus in general, it may be costly to verify this certificate.  But typically the degrees of the covers and the value of $\ell$ will be quite small (much smaller than $g\log q$), in which case computing the $L$-polynomials modulo $\ell$ (or even just enough terms to distinguish them) may be feasible.

\section{Examples}

The simplest case to consider in Conjecture~\ref{main-conj} is when $g=2$ and $q=2$; in this case there are 20 isomorphism classes of curves, all of which have distinct $L$-functions, so one could take $n=0$ in Conjecture~\ref{main-conj}.  The next simplest case is $g=2$ and $q=3$; now there are 69 isomorphism classes of curves, but only 50 isogeny classes of Jacobians.
Of the 50 isogeny classes of Jacobians, 31 contain a unique Jacobian, while 19 contain a pair of Jacobians of non-isomorphic curves.  Among these 19 all but 4 pairs are distinguished by considering the $L$-functions of $H_1(X)$.  These 4 pairs are considered in the first 3 examples below, each of which demonstrates that Conjecture~\ref{main-conj} does not hold if we restrict to $n\le 1$.

\begin{Example}
\normalfont
The genus two curves:
\[
C_1\colon y^2 = 2x^6 + 2x^4 + 2x^3 + 2,\qquad  C_2\colon y^2 = 2x^6 + 2x^5 + x^4 + x^2 + 2x + 2
\]
over $\F_3$ are non-isomorphic, but they have isogenous Jacobians $J_1$, $J_2$ with $L$-polynomial:
\[
9T^4 - 6T^3 + 3T^2 - 2T + 1.
\]

The corresponding Hilbert class fields have degree $\#J_1(\F_3)=\#J_2(\F_3)=5$, and the Riemann-Hurwitz theorem implies that the curves $H_1(C_1)$, $H_1(C_2)$ both have genus 6.
The function fields of $H_1(C_1)$ and $H_1(C_2)$ both have exactly the same number of degree $1,2,3,4,5,6$ places (the counts are $5,0,10,15,60,140$, respectively), which implies that their $L$-polynomials coincide. The computation of
$H_2(C_i)$ seems out of reach so we cannot verify whether these distinguish
the two curves. Instead we look at $2$-power covers in the setting of
Question~\ref{prym-quest}.

The polynomial $f_1(x)$ in the equation $y^2=f_1(x)$ for $C_1$ is irreducible over $\F_3$, while the polynomial $f_2(x)$ in the equation $y^2=f_2(x)$ for $C_2$ splits into irreducible cubic factors; this implies that the Jacobian $J_2$ has full $2$-torsion over $\F_{27}$, while $J_1$ does not. This is already enough to show that the two Jacobians $J_1$ and $J_2$ (and therefore the curves $C_1$ and $C_2$) are non-isomorphic, but this does not immediately fit our approach of computing $L$-polynomials.

However, by taking double covers over $\F_{27}$ and looking at the corresponding elliptic curves, (see Example~\ref{g2} below for a similar calculation), we get elliptic curves with $2$-torsion for the second genus $2$ curve but not for the first, so the isomorphism classes of $C_1$ and $C_2$ are distinguished by the $L$-functions of these double covers. 

\end{Example}

\begin{Example}
\normalfont
The genus two curves:
\[
C_1\colon y^2 = x^5 + x^4 + 2x + 1,\qquad C_2\colon y^2 =  x^5 +x^3 + x^2 + 2x + 2
\]
over $\F_3$ both have $L$-polynomial $L(T) = 9T^4 - 3T^3 + T^2 - T + 1$, and the curves $H_1(C_1)$,$H_1(C_2)$ of genus $8$ also have equal $L$-polynomials.

This example is particularly interesting, in that the corresponding Jacobians $J_1$ and $J_2$ appear to be isomorphic; their respective groups of
$\F_{3^n}$-rational points are isomorphic for $n=1,\ldots,5$. 
As in the previous example, verifying 
Conjecture~\ref{main-conj} seems to be computationally out of reach, but we can distinguish them by taking double covers. We need to work 
over $\F_{3^5}$ and each curve has $15$ \'etale double covers lying in $3$ orbits
of $5$ curves under Frobenius. The Jacobians of the double covers have an
additional elliptic curve factor and we get, as trace of Frobenius for 
these factors, the values (up to sign) of: $28,28,8$ for the first curve and
$28,20,8$ for the second. The appearance of the $20$ shows these elliptic curve factors are 
not isogenous, so the curves are not isomorphic.

A similar example is the pair of curves $C_1\colon y^2 = 2x^6 + x^4 + x^3 + 1$ and $C_2\colon y^2 = x^6 + x^4 + x^3 + 2$ over $\F_3$,
with $L$-polynomial  $L(T) = 9T^4 - 3T^3 + 3T^2 - T + 1$. Traces of Frobenius for 
the elliptic curve factors are  $4,16,28$ for first curve and $4,4,16$ for the second.
\end{Example}

\begin{Example}
\normalfont
The fourth and final example for $g=2$ and $q=3$ is the pair of curves
\[
C_1\colon y^2 = x^5 - 1,\qquad C_2\colon y^2 = x^5 + 1,
\]
which are non-isomorphic quadratic twists.  Their Jacobians are both supersingular with $L$-polynomial $9T^4+1$, and the genus 11 curves $H_1(C_1)$, $H_1(C_2)$ have the same $L$-polynomial.

The curves $C_1$ and $C_2$ both have $4$ points over $\F_3$ and admit a unique (up to twist) unramified double cover. We pin down the double cover by insisting that it have $6$ points over $\F_3$ (the other twist has $2$ points). Then we look at an unramified triple cover of the double cover, of which there are three, all twists of each other. Finally, we see how the 6 points split on these covers and use this information to distinguish the curves.

We have double covers $X_1\colon w^2= x^4+x^3 + x^2 + x +1$ and $X_2\colon  w^2= x^4 - x^3 + x^2 - x +1$ of~$C_1$ and $C_2$ respectively.

Triple covers of $X_1$ are given by $Y_{1,a}: z^3-z = (x+1)w + a$, for $a=0,1,2$.

Triple covers of $X_2$ are given by $Y_{2,a}: z^3-z = (x -1)w + a$, for $a=0,1,2$.

The distinguishing feature is that while all the $Y_{1,a}$ have $\F_3$-points ($12,3,3$, respectively), 
$Y_{2,0}$ is pointless (the curves $Y_{2,1}$ and $Y_{2,2}$ both have $9$ $\F_3$-points).
This implies that the $L$-polynomials of $H_2(C_1), H_2(C_2)$  differ and confirms Conjecture~\ref{main-conj} in this case.
\end{Example}

\begin{Example}
\normalfont
We did an exhaustive search over $\F_2$ and found that there is 
exactly one pair of non-isomorphic smooth plane quartics $C_1,C_2$ over $\F_2$ with the the same $L$-polynomial for which $H_1(C_1),H_1(C_2)$ also have the same $L$-polynomial:
\[
C_1\colon x^3z + xyz^2 + y^4 + y^2z^2 + yz^3,\qquad C_2\colon x^3z + xy^2z + y^4 + y^2z^2 + yz^3.
\]
Both curves have $L$-polynomial $8T^6 - 4T^5 + 2T^3 - T + 1$, with 6 rational points on their Jacobians, and the Hilbert class curves $H_1(C_1)$, $H_1(C_2)$ have genus $13$.

The curves $C_1,C_2$ both have a unique (up to twist) quadratic unramified cover, say $D_1,D_2$ of genus $5$. By the Deuring-Shafarevich formula, $D_1, D_2$ themselves have a unique (up to twist) quadratic unramified cover, and they have distinct $L$-polynomials, even up to quartic twists, which is enough to show $C_1, C_2$ are non-isomorphic and distinguished by the $L$-polynomials of $H_2(C_1)$,$H_2(C_2)$, confirming Conjecture \ref{main-conj} for this example.
\end{Example}

\begin{Example}
\label{g2}
\normalfont
Another example is the pair of genus two curves
\[
C_1\colon y^2 = x^6 +3x^2 + 4,\qquad C_2\colon y^2 = x^6 + 5x^4 + 5x^2 +1
\]
over $\F_7$, which have the same $L$-function. 
To show that they are not isomorphic one can look at the respective double covers and show, by counting points, 
that there cannot be a matching between the double covers of the two curves.
Specifically, both curves have three unramified
double covers defined over $\F_7$.  The Jacobian of these covers
split as the product of the Jacobian of
the original curve with an additional elliptic curve. For the first curve,
all three of these elliptic curves have trace of Frobenius $-4$. For the
second curve, the elliptic curve obtained from the cover $z^2 = x^2+1$
has trace of Frobenius $0$.
\end{Example}

\section{Factoring polynomials over finite fields}

As mentioned in the introduction, the existence of a one-parameter family $X_t$ of curves of genus $g$ over $\F_p(t)$ with
$g$ bounded (or growing very slowly with $p$) such that the $L$-polynomials $L(X_t,T)$ are all distinct (or the number of 
collisions is bounded independent of $p$) for varying $t \in \F_p$ (excluding the $t$ of bad reduction, those for which $X_t$ is singular) leads
to a deterministic polynomial-time algorithm for factoring polynomials in $\F_p[t]$.
There are well-known randomized algorithms to solve this problem whose expected running times are 
polynomially-bounded that are quite fast in practice, so this question is primarily of theoretical interest.
But even for polynomials of degree two, no deterministic polynomial-time algorithm is known, unless one assumes
the Generalized Riemann Hypothesis (GRH), and for general polynomials the question remains open even under GRH.

Using the same heuristic as in section \ref{fund}, there are $p$ choices of $a$ and $p^{g(g+1)/4}$ possible values for the
$L$-polynomial so one would expect this to hold for ``most" families as soon as $g>2$, since $p^2 < p^{g(g+1)/4}$.
Buium \cite{Buium} has shown that most families (in a differential algebraic sense) in characteristic zero 
have finitely many isogeny correspondences, however, even if this result extends to characteristic $p$, it does not rule 
out sporadic isogenies.  Conjecture~\ref{main-conj} does not give the result either, as the genus of the
resulting covers grows too quickly.

We first considered the family $y^2=x^7+(t-1)x^3+tx^2+(t+1)x+1$ of curves of genus~$3$. 
One expects the number of isogeny classes of $3$-dimensional abelian varieties over $\F_p$ to be about $p^3$.
So under our probabilistic heuristic, a one parameter family of curves (with about~$p$ elements) has a probability of about $1/(2p)$ of containing no isogenous Jacobians.
Using the algorithms in \cites{Harvey-Sutherland-I,Harvey-Sutherland-II,Kedlaya-Sutherland} we have verified that for all primes $p\le 10000$ the $L$-polynomials in this family are distinct for all $t$ of good reduction (for each $p$, at most $9$ values of $t \in \F_p$ yield singular curves).
Now $\sum 1/(2p)=O(\log\log p)$ diverges (albeit slowly), so one might expect a collision of $L$-polynomials to occur in this family for some $p > 10^5$ (but one would expect the number of collisions for each $p$ to be bounded by a constant).

To obtain a more compelling example, we instead consider the genus~$4$ hyperelliptic family:
\[
X_t\colon y^2=x^9+(t-1)x^3+tx^2+(t+1)x+1.
\]
Now the number of isogeny classes is on the order of $p^5$, and our heuristic model predicts a probability of roughly $1/(2p^4)$ that two $L$-polynomials $L(X_t,T)$ in our family coincide for some pair of $t\in \F_p$.
The sum $\sum 1/(2p^3)$  now converges.
We have verified that for primes $p\le 2^{17}$ the $L$-polynomials arising in this family are distinct for all~$t$ of good reduction (now at most $11$ values of $t \in \F_p$ yield singular $X_t$), and it seems quite likely that the $L$-polynomials $L(X_t,T)$ arising in the family are distinct for all primes $p$.
Indeed, if  $\pi(t)=t/\log(t)+\varepsilon(t)$ denotes the prime counting function, we can bound the tail of our sum $\sum 1/(2p^3)$ using
\[
\sum_{p>2^{17}} \frac{1}{2p^3} = \int_{2^{17}}^\infty \frac{d\pi(t)}{2t^3} = \int_{2^{17}}^\infty \frac{dt}{2t^3\log t} +  \frac{\varepsilon(t)}{2t^3}\Bigm|_{2^{17}}^\infty  + \int_{2^{17}}^\infty \frac{3\varepsilon(t)dt}{2t^4},
\]
and applying the bound $\varepsilon(t) \le (3t)/(2\log(t)^2)$ (valid for $t\ge 59$) from \cite{Rosser-Schoenfeld} yields
\[
\sum_{p>2^{18}}\frac{1}{2p^3} < 1.187\times 10^{-12} + 4.36\times 10^{-13} + 3.15\times 10^{-13} < 2\times 10^{-12}.
\]
Thus under our heuristic model, the probability that the $L$-polynomials $L(X_t,T)$ at good values of $t$ are not all distinct for every prime $p$ is less than $2\times 10^{-12}$.

These two families were chosen essentially at random by writing a plausible family with no specializations having the
same $L$-polynomial for small primes.
We note that the similar looking families $y^2=x^7+(t-1)x^3+(t+1)x+1, y^2=x^9+(t-1)x^3+(t+1)x+1$ have 
specializations with the same $L$-polynomial for some small primes.

\bigskip

\emph{Acknowledgements:} Both authors would like to
thank Kiran Kedlaya for mentioning the results of \cite{Poonen} at AGCT,
and Bjorn Poonen for writing up that account at our instigation, as well as
Brendan Creutz and Jakob Stix for helpful discussions. The first author 
thanks the National Science Foundation for financial support under grant DMS-152256, and the
second author thanks the Simons Foundation for financial support under grant \#234591 .


\begin{thebibliography}{99}

\bibitem{Buium}
A. Buium,
{\em A finiteness theorem for isogeny correspondences}
in: Journ\'ees de Geometrie Algebrique d'Orsay 1992, Asterisque 218, 1993.

\bibitem{Harvey}
D. Harvey,
{\em Computing zeta functions of arithmetic schemes},
Proc. Lond. Math. Soc. {\bf 111} (2015) 1379--1401.

\bibitem{Harvey-Sutherland-I}
D. Harvey and A.V. Sutherland
{\em Computing Hasse-Witt matrices of hyperelliptic curves in average polynomial time}, LMS J. Comput. Math. {\bf 17} (2014) 257--273.

\bibitem{Harvey-Sutherland-II}
D. Harvey and A.V. Sutherland
{\em Computing Hasse-Witt matrices of hyperelliptic curves in average polynomial time.}
In: D. Kohel and I. Shparlinski (Eds.): Frobenius distributions: Lang-Trotter and Sato-Tate conjectures, Contemp. Math. 663, American Mathematical Society, 2016.

\bibitem{Hess}
F. Hess, 
{\em An algorithm for computing isomorphisms of algebraic function fields.} 
In: Buell, D. (Herausgeber): Proceedings of the Sixth International Algorithmic Number Theory Symposium, 
ANTS-VI, LNCS 3076, Springer-Verlag, Berlin-Heidelberg-New York. 2004.

\bibitem{Howe-DiPippo}
E. Howe and S. DiPippo
{\em Real polynomials with all roots on the unit circle and abelian varieties over finite fields}
J. Number Theory {\bf 73} (1998) 426--450; Corrigendum, J. Number Theory {\bf 83} (2000) 182.

\bibitem{Kedlaya-Sutherland}
K. Kedlaya and A.V. Sutherland
{\em Computing $L$-series of hyperelliptic curves.}
In: A.J. van der Poorten and A. Stein (Eds.): Proceedings of the Eighth International Algorithmic Number Theory Symposium,
ANTS-VIII, LNCS 5011, Springer-Verlag, Berlie-Heidelberg, 2008.

\bibitem{Stirpe}
P. Mercuri and C. Stirpe, 
{\em Classification of algebraic function fields with class number one},
J. Number Theory, {\bf 154}, (2015), {365--374}.

\bibitem{Mestre}
J.-F. Mestre,
{\em Couples de jacobiennes isog\`enes de courbes hyperelliptiques de genre arbitraire}
arxiv 0902.3470

\bibitem{Moehlmann}
G. M\"ohlmann,
{\em Einbettungen globaler Funktionenk\"orper}
Diplomarbeit TU Berlin 2008

\bibitem{Poonen}
B. Poonen,
{\em Using zeta functions to factor polynomials over finite
fields}, Preprint, 2017.

\bibitem{Rosser-Schoenfeld}
J. Rosser and L. Schoenfeld, {\em Approximate formulas for some functions of prime numbers}, Illinois J. Math. \textbf{6} (1962), 64--94.

\bibitem{Shankar-Tsimerman}
A. S. Shankar and J. Tsimerman,
{\em Unlikely intersections in finite characteristic}
arxiv 1610.03552

\bibitem{Shen}
Q. Shen and S. Shi, 
{\em Function fields of class number one}, 
J. Number Theory {\bf 154} (2015) 375 -- 379.
 
\bibitem{Tuitman}
J. Tuitman,
{\em Counting points on curves using a map to $\mathbb{P}^2$, II},
Finite Fields Appl. {\bf 45} (2017) 301--322.

 
\end{thebibliography}
\end{document}